\documentclass[11pt]{amsart}
\usepackage{amsmath,amsfonts}
\usepackage{amssymb}
\usepackage{amscd}
\usepackage{amsthm}
\usepackage{yhmath}
\usepackage{subfigure}

\usepackage[all]{xy}
\usepackage{color}

\usepackage{marginnote}

\numberwithin{equation}{section}

\newtheorem{theorem}[equation]{Theorem}

\newtheorem{proposition}[equation]{Proposition}

\theoremstyle{remark}
\newtheorem{remark}[equation]{Remark}

\theoremstyle{definition}

\def\XXint#1#2#3{{\setbox0=\hbox{$#1{#2#3}{\int}$}
	\vcenter{\hbox{$#2#3$}}\kern-.5\wd0}}

\newcommand{\N}{\mathbb N}

\newcommand{\R}{\mathbb R}

\def\eps{\epsilon}
\def\ga{\gamma}

\def\la{\lambda}

\begin{document}

\title[Non-minimality of corners]{
Corners in non-equiregular sub-Riemannian manifolds} 

\author[Le Donne]{Enrico Le Donne}

\address[Le Donne]{Department of Mathematics and Statistics, P.O. Box 35,
FI-40014,
University of Jyv\"askyl\"a, Finland}%
\email{ledonne@msri.org}

\author[Leonardi]{Gian Paolo
Leonardi}
\address[Leonardi]{Universit\`a di Modena e Reggio Emilia,
  Dipartimento di Scienze
Fisiche, Informatiche e Matematiche, via Campi 213/b, 41100
  Modena, Italy}

\email{gianpaolo.leonardi@unimore.it}

\author[Monti]{Roberto Monti}
\address[Monti and Vittone]{Universit\`a di Padova, Dipartimento di Matematica,
via Trieste 63, 35121 Padova, Italy}
\email{monti@math.unipd.it}

\author[Vittone]{Davide Vittone}
\email{vittone@math.unipd.it}


\renewcommand{\subjclassname}{%
 \textup{2010} Mathematics Subject Classification}

 \date{20-2-2104}

\begin{abstract} 
We prove that in a class
of non-equiregular sub-Riemannian manifolds 
corners are not length minimizing. This extends the results of
\cite{Leonardi-Monti}.
As an application of our main result we complete and simplify the analysis in \cite{Monti-nonminimizing}, showing that in a
$4$-dimensional  sub-Riemannian structure suggested by Agrachev and Gauthier all length-minimizing curves are smooth.
 \end{abstract}

\maketitle


 

\section{Introduction}
\setcounter{section}{1}

One of the major open problems in sub-Riemannian geometry 
is the regularity of
length-minimizing curves.
Indeed, no example  of a non-smooth minimizer is known, 
and even the possibility
of minimizers with singularities of corner-type 
has not yet been excluded in full generality 
(see  Problem II in \cite{Agrachev_problems} and the
discussion in Section 4 of
\cite{Monti_survey}).

In \cite{Leonardi-Monti}, the second and third-named 
authors introduced a
shortening technique specifically designed for showing 
the non-minimality of 
curves with corner-type singularities (see also the developments in 
\cite{Monti_CV}). 
This technique works for a class of equiregular
sub-Riemannian manifolds satisfying the technical condition
\eqref{LLMM} below.

In this paper, we prove the non-minimality 
of corners in a class of sub-Riemannian
manifolds of non-equiregular type. Namely, we show
that if the horizontal
distribution satisfies condition \eqref{LALLA} below at the
corner point, then the curve is not length
minimizing.
In this case, the construction of a competitor shorter than the corner
is simpler than the one in \cite{Leonardi-Monti}
and relies on the Nagel-Stein-Wainger estimates \cite{nagelstwe} for the
Carnot-Carath\`eodory distance.

Let $M$ be an $n$-dimensional differentiable manifold
and let $\mathcal D\subset TM$ be an $m$-dimensional distribution of planes,
for some $2 \leq m\leq n$.
Then, $\mathcal D(x)\subset T_xM$ is 
an $m$-dimensional subspace of the tangent space $T_xM$,
for all $x\in M$. 
Let $X_1,\ldots, X_m$ be a frame of vector fields
that form a basis
for $\mathcal D(x)$, that is 
$\mathcal D(x) = \mathrm{span} \{ X_1(x),\ldots, X_m(x)\}$
for each $x\in M$.
This frame always exists locally.
The assumption that the fiber $\mathcal D(x)$
has a constant dimension on $M$ 
plays no role in our argument and can be dropped.

We denote by $\mathcal I = 
\bigcup _{i \geq 1} \{ 1,\ldots,m\}^ i 
$ the set of admissible multi-indices.
For any  $\beta = (\beta_1,\ldots,\beta_i)\in \mathcal I$, 
for some $i\geq 1$, 
let us define the iterated commutator 
\begin{equation}\label{pax2}
     X_\beta = [X_{\beta_i},[X_{\beta_{i-1} },\ldots 
[X_{\beta_2},X_{\beta_1}]\ldots]].
\end{equation}
We say that $\mathrm{length}(\beta)=i$ is the length of the multi-index
$\beta$.
Analogously, we say that $
\mathrm{length}(X_\beta)=i$
is the length of the iterated commutator $X_\beta$.
For any point $x\in M$ and $i\geq 1$, let
\[
\mathcal D_i(x) = \mathrm{span}\big\{ X_\beta(x)\in T_xM \, : \,  
\mathrm{length}(X_\beta) =i\big\}.
\] 
Finally, we let $\mathcal L _i(x) = 
\mathcal D_1(x)+\ldots+\mathcal D_i(x)$ for $i\geq 1$, and we also
agree that  $\mathcal L_0(x) = \{ 0\}$.
We assume that $\mathcal D$ is bracket generating, i.e., 
for any $x\in M$ there exists an index $i\in\N$ 
such that $\mathcal L_i(x)= T_x M$.

An absolutely continuous curve $\gamma:[0,1]
\to M$ is said to be
\textit{horizontal}  with respect to the distribution $\mathcal D$ (or simply $\mathcal D$-horizontal) if  there exist  bounded measurable 
functions $h_{1},\ldots,h_{m}:[0,1]\to \R$ such
that
\[
\dot\gamma(t) = \sum_{j=1}^ m h_{j}(t)\, X_{j}(\gamma(t)),\quad
\text{for almost every } t \in [0,1].
\]
Let 
$g(x;\cdot)$ be a positive quadratic form (metric) 
on $\mathcal D(x)$, $x\in M$. 
The  length of $\gamma$ in the sub-Riemannian manifold
$(M,\mathcal D, g)$ is defined as
\[
     L(\gamma) = \int_0^1\sqrt{ g( {\gamma(t)} ; \dot\gamma(t)) } \, dt,
\]
and the {\em sub-Riemannian distance} 
between two points  ${x}, {y}\in M$
is defined as 
\[
d (x, y) := \inf\big\{ L(\gamma):\ \gamma\in AC([0,1];M) \text{
horizontal}, \ \gamma(0)= x,\  \gamma(1)=y\big\}.
\]
When $M$ is connected, the above set  is always nonempty because the 
distribution $\mathcal D$ is bracket-generating, 
and $d$ is a distance on
$M$.
Finally, we say that a horizontal 
curve $\gamma$ joining $x$ to $y$ minimizes the sub-Riemannian
length (i.e., it is a length minimizer) if $L(\gamma) = d(x,y)$.

Let $\gamma:[0,1]\to M$ be a $\mathcal D$-horizontal curve.
When they exist, we  denote by $\dot\gamma_L(t)$ and
$\dot\gamma_R(t)$ the left and right derivative
of $\gamma$ at the point $t\in(0,1)$. 
We say that  $\gamma$ has
a corner at
the point $x=\gamma(t)\in M $, if the left and right derivatives
at $t$,
 do exist and are linearly independent. 
In \cite{Leonardi-Monti}, it is shown
that if the distribution $\mathcal D$ is equiregular (i.e., for every $i\geq 1$ the dimension of $\mathcal D_{i}(x)$ is constant on $M$) and satisfies the condition
\begin{equation}\label{LLMM}
   [\mathcal D_i,\mathcal D_j]
  \subset  \mathcal L_{i+j-1},\quad \textrm{for $i,j\geq2$ such that $i+j>4$},
\end{equation}
then corners in $(M,\mathcal D, g)$ are not length minimizing. 
In this paper, we prove that if the distribution $\mathcal D$ satisfies at some point $x\in M$ the condition
\begin{equation}
\label{LALLA}
   \mathcal L_i(x)\neq  \mathcal L_{i-1} (x) 
\quad \Rightarrow\quad \mathcal L _{i+1} (x) = \mathcal L_i(x),\qquad \text{for all $i\geq 2$},
\end{equation}
then corners at $x$ are not length minimizing.

\begin{theorem}
 \label{thermos2} 
Let $\gamma:[0,1]\to M$ be a horizontal
curve with a 
corner at the point $x=\gamma(t)\in M$, for  
$t\in(0,1)$. If the distribution $\mathcal D$ satisfies \eqref{LALLA}
at  $x$, then  
$\gamma$ is not length minimizing 
in $(M,\mathcal D, g)$.
\end{theorem}

The proof of Theorem \ref{thermos2}, the main result of this paper, is presented in Section
\ref{S2}.
After a blow-up
argument, we can assume that $M=\R^n$, that $\mathcal D$ is
a $2$-dimensional distribution of planes in $\R^n$,
and that    $\gamma:[-1,1]\to\R^n$ is a  
corner at the point $0\in\R^n$
of the type
\[
\gamma(t) = \begin{cases}
-t x  & \text{if }t\in [-1,0]\\ 
\,\,\,\,t y & \text{if }t\in (0,1],
\end{cases}
\] 
where $x,y\in \R^n$ are linearly independent.
We prove the non-minimality of $\gamma$ by an 
inductive argument on the dimension $n\geq 2$.
In the inductive step, we use assumption
\eqref{LALLA} and known estimates on the sub-Riemannian distance
to find a competitor shorter than the corner.

We found the basic idea of the proof of Theorem \ref{thermos2}
starting from a question
raised by A.~Agrachev and J.~P.~Gauthier during the meeting 
\emph{Geometric control and sub-Riemannian geometry} held in
Cortona in May 2012.
They suggested  the
following situation, in order to find a   \emph{nonsmooth}
length-minimizing curve. On the manifold $M=\R^4$, let $\Delta$ be the  
distribution of 2-planes spanned point-wise by the 
vector fields 
\begin{equation}
 \label{X_1,X_2}
     X_{1} = \frac{\partial }{\partial x_1}  +2x_{2} \frac{\partial }{\partial
x_3}  + x_{3}^{2} \frac{\partial }{\partial x_4} ,\qquad 
     X_{2} =  \frac{\partial }{\partial x_2}  - 2x_{1} \frac{\partial }{\partial
x_3} .
\end{equation}
%
%
The distribution $\Delta$ satisfies \eqref{LALLA}. 
We fix on $\Delta$ the quadratic form $g$ making $X_1$ and $X_2$ orthonormal.

Let  $\alpha>0$ be a parameter and 
consider the initial and final points $y =
(-1,\alpha,0,0)\in\R^4$ and  $x = (1,\alpha,0,0)\in\R^4$, respectively.
Agrachev and Gauthier asked whether the corner $\gamma:[-1,1]\to\R^4$
joining $y$ to $x$ 
\begin{equation}\label{barba}
   \gamma_1(t) = t, \quad
  \gamma_2(t) = \alpha|t|,\quad  
\gamma_3(t) = 0,\quad
 \gamma_4(t)=0, \quad t\in[-1,1]
\end{equation}
is, for small $\alpha>0$, a length minimizer in $(\R^4,\Delta, g)$.
The presence of the variable $x_{3}$ 
in the coefficients of the vector field
$X_{1}$ in \eqref{X_1,X_2}  is the 
technical obstruction for the application of
the results  of \cite{Leonardi-Monti}.

In \cite{Monti-nonminimizing}, the curve $\gamma$ in
\eqref{barba} was shown not to be length minimizing for $\alpha\neq 1$,
by the explicit construction of a shorter 
competitor. This answered the above question in the
negative.
The case $\alpha =1$, however, was left 
open.

In Section
\ref{S3}, as an application of Theorem \ref{thermos2}, we prove the following result.

 \begin{theorem}\label{main teo} Let $g$ be any smooth metric on $\Delta$. 
In the sub-Riemannian manifold
$(\R^4,\Delta,g)$
all length minimizing curves are smooth and, in particular,
no corner is length minimizing. 
\end{theorem}

The proof of Theorem \ref{main teo} relies on Theorem \ref{thermos2}.
The inductive base
is provided by the regularity 
of geodesics in the first Heisenberg group. 
This proof covers in particular the case $\alpha=1$ in \eqref{barba} 
and is simpler than the one in \cite{Monti-nonminimizing}.

\section{Proof of Theorem \ref{thermos2}}
\label{S2}

The first step of the proof is a blow-up argument that
closely follows \cite{Leonardi-Monti}.

Let $\gamma:[-1,1]\to M$ be a horizontal curve with a corner at the point
$x=\gamma(0)\in M$. 
We can choose smooth and linearly independent vector fields $X_1,X_2\in\mathcal D$ such
that $X_1(x) = \dot\gamma_R(0)$ and $X_2(x) = -\dot\gamma_L(0)$
and we complete
$X_1,X_2$ to a (local) frame $X_1,\ldots,X_m$ for $\mathcal D$.
Then we complete $X_1,\ldots,X_m$ 
to a frame  $X_1,\ldots, X_n$ for $TM$ in the following way.
We choose iterated commutators $X_{m+1},\ldots,X_n \in \{ X_\beta\,:\, 
 \beta \in \mathcal I,\, \textrm{length}(\beta)\geq 2 \}$
such that $X_1(x),\ldots, X_n(x)$ are linearly independent.
This choice is possible because $\mathcal D$ is bracket generating at $x$.
We can also assume that $j\leq k$ implies $\mathrm{length}(X_j)\leq
\mathrm{length}(X_k)$.

In a neighbourhood of $x\in M$, we fix exponential coordinates of the
first type induced by the frame $X_1,\ldots,X_n$ 
starting from $x$. Then
we can identify $M$ with $\R^n$, $X_1,\ldots,X_n$ 
with vector fields on $\R^n$,
and $x$ with $0\in\R^n$. 
The fact that we have exponential coordinates
of the first type means that for 
$ x = (x_1,\ldots,x_n)\in \R^n $ (in fact, for
$x$ belonging to a neighbourhood of $0\in \R^n$) we have 
\begin{equation}
 \label{poll}
 x = \exp\Big(\sum_{i=1}^n x_i X_i\Big)(0) .
\end{equation}
Here, the exponential mapping is defined by $\exp(X)(0) = \gamma(1)$ where
$\gamma$ is the solution of $\dot\gamma = X(\gamma)$ and $\gamma(0)=0$.

We assign to the coordinate $x_i$ 
the weight $w_i=\mathrm{length}(X_i)$, $i=1,\ldots,n$.
Then we have $w_1=\ldots= w_m=1$. The natural dilations on $\R^n$ adapted to the
frame $X_1,\ldots,X_n$ are
\begin{equation}\label{dilations2}
 \delta_\lambda(x) = (\la ^{w_1} x_1, \la^{w_2}  x_2, \ldots,   \la^{w_n} 
x_n),\quad
x\in\R^n,\,\, \la>0.
\end{equation}
Let $X=X_\beta$ be any iterated commutator of the vector fields 
$X_1,\ldots,X_m$. Then we have
\begin{equation}\label{forx}
  X = \sum_{i=1}^n a_i (x) \frac{\partial}{\partial x_i},
\end{equation}
where $a_i\in C^\infty(\R^n)$, $i=1,\ldots,n$, are smooth functions that
have the structure described in the following proposition.

\begin{proposition}\label{ciox}
There exist polynomials $p_i:\R^n\to\R$ and
functions $r_i :\R^n\to\R$, $i=1,\ldots,n$, such that:
\begin{itemize}
 \item[i)] $a_i(x) = p_i(x) + r_i(x)$, $x\in\R^n$;

\item[ii)] $p_i(\delta_\lambda(x)) = \lambda ^{w_i-\mathrm{length}(X)} p_i(x)$;

\item[iii)] $\displaystyle \lim_{\la\to\infty} \lambda ^{w_i
-\mathrm{length}(X)}
r_i(\delta _{1/\lambda} (x))=0$, $x\in\R^n$.
\end{itemize}
\end{proposition}

\noindent Proposition \ref{ciox}
can be proved as in \cite[page 306]{Margulis-Mostow2}.
We omit the details, here. For $\lambda>0$, we let
\begin{equation}\label{NUMO}
 X^\lambda (x) = \sum_{i=1}^n \lambda^{w_i -\mathrm{length}(X) } a_i
(\delta_{1/\lambda} (x))\frac{\partial}{\partial x_i},\quad x\in\R^n.
\end{equation}
The mapping $X\mapsto X^\lambda$ is bracket-preserving.
Namely, for any multi-index $\beta\in\mathcal I$ and for $i=1,\ldots,m$ we have
\begin{equation}\label{BG}
 [X_i,X_\beta]^\lambda = [X_i^\lambda, X_\beta^\lambda],\quad \lambda>0.
\end{equation}
We let $\mathcal D^\lambda 
= \mathrm{span}\{ X_1^\lambda,\ldots,X_m^\lambda\}$,
$\mathcal D^{\lambda}_i = \mathrm{span} \{ X_\beta ^\lambda\, 
: \, \mathrm{length}(\beta)=i\big\}$,
and $\mathcal L^{\lambda}_i =\mathcal D_1^\lambda +\ldots+ \mathcal D_i^\lambda$.
By \eqref{BG}, from \eqref{LALLA} we deduce that at the point $x=0$ we have
\begin{equation}
\label{LALLA2}
   \mathcal L_i^\lambda \neq  \mathcal L_{i-1} ^\lambda 
\quad \Rightarrow\quad \mathcal L _{i+1} ^\lambda = \mathcal L_i^\lambda,
\quad
  \textrm{for $i\geq 2$}.
\end{equation}

By Proposition \ref{ciox}, for any iterated commutator $X=X_\beta$ as in
\eqref{forx},
we can define the vector field $X^\infty$ in $\R^n$
\[
 X^\infty(x) = \lim_{\la\to\infty}  X^\lambda (x) = \sum_{i=1}^n p_i(x)
\frac{\partial}{\partial x_i},\quad x\in\R^n,
\]
where $p_i$, $i=1,\ldots,n$, are polynomials such that
$p_i\circ\delta_\lambda
= \lambda ^{w_i-\mathrm{length}(X)}p_i$. In particular, if $w_i<\mathrm{length}
(X)$ then $p_i = 0$. Passing to the limit as $\lambda\to\infty$ in \eqref{BG},
we see that also the mapping $X\mapsto X^\infty$ 
is bracket-preserving, i.e., $
 [X_i,X_\beta]^\infty = [X_i^\infty, X_\beta^\infty]$.
Then at the point $x=0$, 
condition \eqref{LALLA2} holds also for $\lambda=\infty$.

Let $g(x;\cdot)$ be   a metric   on $\mathcal D(x)$.
On the distribution $\mathcal D^\la$, $\la>0$, 
we introduce the metric $g^\la(x;\cdot)$
defined by 
\[
  g^\lambda (x;X^\lambda ) = g(\delta_{1/\lambda} (x); X),\quad
x\in\R^n,
\]
and  on
$\mathcal D^\infty = \mathrm{span} \{X_1^\infty,\ldots ,X_m^\infty\}$ we
introduce the metric $g^\infty(x;\cdot)$ defined by
\[
          g^\infty(x; X^\infty ) =
          \lim_{\lambda \to \infty } g^\lambda (x;
          X^\lambda) = g(0;X),\quad x\in\R^n.
\]

We blow up the curve $\ga$ at the corner   point $0\in\R^n$. 
For $\la>0$ and $t
\in[-\la,\la]$, let  $\ga^\la(t) = \delta _{\la} \ga( t/\la)$.  
Because $\dot\gamma_R(0) = X_1(0)$, $\dot\gamma_L(0)=-X_2(0)$, 
we obtain the limit curve $\displaystyle \gamma^\infty = \lim_{n\to\infty}
\gamma^\lambda$, 
\begin{equation}
 \label{colle}
       \ga^\infty(t) =
    \left\{
  \begin{array}{rl}
  \displaystyle
    \mathrm{e}_1 t, & t\in[0,1] \\
  \displaystyle
     -\mathrm{e}_2 t, & t\in [-1,0),
  \end{array}
\right.
\end{equation}
where $\mathrm{e}_1 = (1,0,\ldots,0)$ and $\mathrm{e}_2=(0,1,0,\ldots,0)$.

\begin{proposition}\label{MINO}
If the curve $\gamma$ is length minimizing in $(M,\mathcal D, g)$
then the curve $\gamma^\infty$ is length minimizing
in $(\R^n,\mathcal D^\infty,g^\infty)$.
\end{proposition}

Proposition \ref{MINO} is proved in \cite{Leonardi-Monti}, Proposition 2.4.
Our goal is to prove that the corner $\gamma^\infty$ is not length minimizing
in $(\R^n,\mathcal D^\infty, g^\infty)$.
Thus, we can without loss of generality assume that $M=\R^n$, $\mathcal D
= \mathcal D^\infty$, and $\gamma=\gamma^\infty$.
Since   $\gamma$ is contained in the orbit of the distribution
$ \mathrm{span}\{ X_1, X_2\}$,
we can also  assume that $m=2$.
Finally, we can pass to exponential coordinates of the 
second type associated with $X_1,\ldots,X_n$. Namely, we can assume that
for all $x = (x_1,\ldots,x_n) \in \R^n$ belonging to a neighbourhood of $0\in\R^n$ we have
\[
     x = \exp(x_1 X_1) \circ \ldots \circ \exp (x_n X_n) (0).
\]
Then we can also assume that $X_1$ and $X_2$ are vector fields in $\R^n$ of the form  
\begin{equation}\label{X!}
\begin{split}
  X_{1}  =  \frac{\partial }{\partial x_1} 
\quad\textrm{and}\quad
  X_{2} =  \frac{\partial }{\partial x_2}  + \sum_{i=3}^{n} p_{i}(x)
 \frac{\partial }{\partial x_i} ,
\end{split}
\end{equation}
where $p_{i}:\R^n\to\R$, $i=3,\ldots,n$, are  polynomials 
of the variable $x = (x_{1},\ldots,x_{n})\in\R^n$ such that 
$p_i(\delta_\lambda(x)) = \lambda ^{w_i-1} p_i(x)$.

Condition \eqref{LALLA2} passes to the limit as $\lambda\to\infty$.
Then,   assumption  \eqref{LALLA} at
the point $x=0$ reads
\begin{equation}
\label{LALLAbis}
   \mathcal L_i(0)\neq  \mathcal L_{i-1} (0) 
\quad \Rightarrow\quad \mathcal L _{i+1} (0) = \mathcal L_i(0),\quad
  \textrm{for all $i\geq 2$}.
\end{equation}
Condition   $\mathcal L _{i+1}(0) = \mathcal L  _i(0) $ is equivalent to
$\mathcal D_{i+1}(0)\subset \mathcal L_i(0)$. 
As $\mathcal D=\mathcal D^\infty$,
\eqref{LALLAbis} is
equivalent to
\begin{equation}
\label{DELTA}
   \mathcal L_i(0)\neq  \mathcal L_{i-1} (0) 
\quad \Rightarrow\quad   \mathcal D _{i+1} (0) = \{0\},\quad
  \textrm{for all $i\geq 2$}. 
\end{equation}
We prove this claim. Let $X_\beta$ be an iterated commutator such that
$\mathrm{length}(\beta) =
i+1$, i.e., $X_\beta\in\mathcal D_{i+1}$. According to  Proposition \ref{ciox},
we have
\[
  X_\beta  = \sum_{j=1}^n p_j (x) \frac{\partial}{\partial x_j},
\]
where $p_j$ are polynomials 
satisfying  $p_j(\delta_\lambda(x)) = \lambda ^{w_j-i-1 }p_j(x)$.
Then the sum above ranges over indices $j$ such that $w_j\geq i+1$.
On the other hand, $\mathcal L_{i}(0) =
\mathrm{span}\big\{
\frac{\partial}{\partial x _j} : w_j \leq i \big\} $ and thus if  $X_\beta(0) \in \mathcal
L_i(0)$ we conclude that $X_\beta(0) = 0$. This proves \eqref{DELTA}.

For given indices  $j,k=1,\ldots,n$, 
we say  that 
$j\prec k$ if there exists $i\geq 2$ such that
 $X_j(0) \in \mathcal D_i(0) $ but   $X_k(0)\notin \mathcal L _i(0)$.
From \eqref{DELTA}, we deduce that the weights $w_1,\ldots,w_n$
satisfy
the following condition:
\begin{equation}
\label{lenghe}
j \prec k\quad \Rightarrow\quad   
       w_{j}+2\leq w_k.
\end{equation}

We are ready to prove that  the corner $\gamma=\gamma^\infty$
in \eqref{colle}
 is not a length minimizer
in $(\R^n,\mathcal D, g)$; this will give a contradiction and prove Theorem
\ref{thermos2}. 

We can without loss of generality assume that $X_1$ and $X_2$ 
in \eqref{X!} are orthonormal with respect to the metric $g$. 
This is because two different metrics are locally equivalent
and the equivalence constants do not affect our estimates below. 
Then the length of a $\mathcal D$-horizontal curve $\gamma:[0,1]\to\R^n$ is
\begin{equation}
  L(\gamma) = \int_0^1\sqrt{  \dot\gamma_{1}(t)^2+ \dot\gamma_{2}(t)^2}  dt.
\end{equation}

The proof is by induction on the dimension $n$ of $\R^n$.
In order to fix the base of induction we distinguish two cases:

\begin{itemize}
 \item[1)] We have $\mathcal L_2(0) = \mathcal L_1(0)$. In this case,
the base of induction is $n=2$. On $\R^2$ we have the standard Euclidean metric
and corners are not length minimizing.

 \item[2)] We have $\mathcal L_2(0) \neq \mathcal L_1(0)$.
In this case,
the base of induction is $n=3$. On $\R^3$ we have the 
Heisenberg group structure. We know that corners are not length minimizing
for any sub-Riemannian metric in the Heisenberg group.
\end{itemize}


We assume that the claim holds for $n-1$
with $n\geq 3,4$ in the two cases,
and   we  prove it for $n$.

Let $\pi:\R^n\to\R^{n-1}$ 
be the  projection $\pi(x_1,\ldots,x_n) = (x_1,\ldots,x_{n-1})$,
and define the vector fields $\widehat X_1 = \pi_* X_1$ and 
$\widehat X_2 = \pi_* X_2$.
Recall that for each $j=3,\ldots,n$ the polynomial $p_j$
appearing in $X_2$ in \eqref{X!}
satisfies $p_j\circ\delta_\lambda = \lambda^{w_j - 1} p_j$ and thus
it depends only on the variables $x_1,\ldots, x_{j-1}$.
In particular, for each $j=3,\ldots,n$ the polynomial
$p_j$  does not depend on $x_n$. 
It follows that
\[
\begin{split}
\widehat   X_{1}(x)  =   \frac{\partial}{\partial x_1} ,
\quad
\widehat   X_{2}(x)  =   \frac{\partial}{\partial x_2} + \sum_{j=3}^{n-1} 
    p_{j}(x)  \frac{\partial}{\partial x_j},\quad \textrm{where } x =
(x_1,\ldots,x_{n-1}) \in\R^{n-1}.
\end{split}
\]
We let  
 $\widehat{ \mathcal{ D
}} = \mathrm{span}
\{ \widehat X_1,\widehat X_2\}$
and we denote by $\widehat g$ 
the metric on $\widehat{ \mathcal D}$ that makes $\widehat X_1$ and $\widehat X_2$ orthonormal.
The distribution $\widehat {\mathcal D} $  satisfies \eqref{LALLA}.

The projection of the curve  $\gamma$ in \eqref{colle}
to $\R^{n-1}$, the curve $\widehat\gamma = \pi(\gamma) =
(\gamma_1,\dots,\gamma_{n-1})$,
is a corner at $0\in\R^{n-1}$. 
By the inductive assumption, this curve is not length minimizing 
in $(\R^{n-1}, \widehat {\mathcal D}, \widehat g)$.
Then there exists a $\widehat{ \mathcal D}$-horizontal 
curve $\widehat\sigma =
(\widehat \sigma_1, \ldots,\widehat \sigma_{n-1})$ in $\R^{n-1}$ joining 
the point   $\mathrm{e}_2\in \R^{n-1}$ to the point  $\mathrm{e}_1\in
\R^{n-1}$ 
and  satisfying  
\[
        k := L(\widehat\gamma)-L(\widehat\sigma) >0.
\]

Let  $\sigma:=(\sigma_1, \ldots,\sigma_n) $ 
be the $\mathcal D$-horizontal lift to $\R^n$ 
of the plane curve
$(\widehat \sigma_{1},\widehat \sigma_{2})$ 
starting from the initial point   $\mathrm{e}_2$.
Clearly, we have $\sigma_i =\widehat \sigma_i$ for $i=1,\ldots,n-1$
and
\[
              L(\sigma ) 
         = \int_0^1\sqrt{  \dot\sigma_{1}(t)^2+ \dot\sigma_{2}(t)^2}  dt
            = L(\widehat \sigma ).
\]
Finally, the end-point of $\sigma$ is of the form 
$\mathrm{e}_1+ h \mathrm{e}_n \in\R^n$, for
some $h\in \R$.

By our choice of the basis $X_1,\ldots,X_n$, 
there exists a multi-index
$\beta= (\beta_1,\ldots,\beta_i)\in\mathcal I $, $i\geq 3$, such that
$X_n= X_\beta$. Since we are in exponential coordinates and also using
Proposition
\ref{ciox}, we deduce that $X_n=\partial /\partial x_n$.
Thus, we have
\begin{equation}\label{pax}
      \frac{\partial}{\partial x_n} = X_\beta =
[X_{\beta_i},[X_{\beta_{i-1} },\ldots 
[X_{\beta_2},X_{\beta_1}]\ldots]].
\end{equation}
The integer $w_n = i $ is the length of the multi-index.
We define  the multi-index
 $\widehat \beta  =  (\beta_1,\ldots,\beta_{i-1})$, that has  
length $i-1 = w_n-1$, and we define the corresponding iterated commutator
\[ 
Z = X_{\widehat \beta} = 
[X_{\beta_{i-1} },\ldots 
[X_{\beta_2},X_{\beta_1}]\ldots]=\sum_{j=1}^n b_j(x)   \frac{\partial}{\partial
x_j},
\]
where $b_j\in C^\infty(\R^n)$ are suitable functions, and, in fact,
polynomials.
By Proposition \ref{ciox}, these polynomials are homogeneous:
\[
 b_j (\delta_\lambda(x)) 
= \lambda ^{w_j-w_n+1} b_j(x),\quad x\in\R^n.
\]
%
%
%
%
%
%
Thus, when $w_j-w_n+1 <0$   the polynomial
$b_j$ vanishes identically, $b _j=0$, and  
the vector field $Z$ has the form
\[
  Z =  \sum_{w_j \geq w_n-1} b_j (x)   \frac{\partial}{\partial x_j}.
\]
If $w_j=w_n-1$ then $b_j(x)$ has homogeneous degree $0$ and thus
it is constant. On the other hand, 
we have $ {\partial} / {\partial x_n}\in\mathcal D_i(0)$ and thus
$\mathcal D_i(0)\neq \{0\}$. From \eqref{DELTA} it follows that
$\mathcal L_{i-1}(0) =\mathcal L_{i-2}(0)$, 
that is $\mathcal D_{i-1}(0) = \{0\}$.
Because we have $Z\in \mathcal D_{i-1}$, then $Z(0)=0$ and we conclude that
$b_j=0$ when $w_j = w_n-1$ and 
$Z$ is, in fact, of the form
\[
  Z =  \sum_{w_j = w_n } b_j (x)   \frac{\partial}{\partial x_j},
\]
with $b_{j}(\delta_{\lambda}(x)) = \lambda b_{j}(x)$. Therefore we have 
\begin{equation}\label{formabj}
b_{j}(x) = c_{j1}\,x_{1}+c_{j2}\,x_{2}
\end{equation}
 for all $j$ such that $w_{j}=w_{n}$, and for suitable constants $c_{j1},c_{j2}$. 
Since the coefficients of $X_2$ (and $X_1$) in \eqref{X!} 
do not contain the
variables $x_j$ 
such that $w_j= w_n$, we infer that
\[
    \frac{\partial}{\partial x_n}= [X_{\beta_i}, Z] 
= \sum_{w_j=w_n} \partial _{\beta_i} b_j(x)  \frac{\partial}{\partial x_j},
\]
and this implies that 
\begin{equation}\label{bubbles}
   \frac{\partial}{\partial x_{\beta_i} } b_n(x) = 1
     \qquad \textrm{and} \qquad
         \frac{\partial}{\partial x_{\beta_i} } b_j(x) = 0,\quad j\neq n,
\end{equation}
 where
either $\beta_i=1$ or $\beta_i=2$. We conclude that either $c_{n1}=1$ or $c_{n2}=1$ (or both).

Assume that $c_{n1}=1$. The proof in the case $c_{n2}=1$ is analogous. 
By our choice of the basis $X_1,\ldots,X_n$,
for any $j=3,\ldots, n-1$ there exists
a multi-index $\beta^j\in\mathcal I$ such that
\[  
 X_j = X_{\beta^{j}}  =  \frac{\partial}{\partial x_j}
+\sum_{k=j+1}^n p_{jk}(x)   \frac{\partial}{\partial x_k},
\]
for suitable polynomials $p_{jk}$. Thus, at the point $x= \mathrm{e}_1\in\R^n$,
the vectors
\[
X_1(x), X_2(x),\ldots ,X_{n-1}(x), Z(x)
\]
are linearly independent, i.e., they  form a basis
of $T_x \R^n$. In particular, the vector field $Z$
is an iterated commutator of $X_1$ and $X_2$ with length $w_n-1$.
By the Nagel-Stein-Wainger estimate 
for the Carnot-Carath\`eodory distance
(see  \cite{nagelstwe} and, in particular, Theorem 4),
 there exist a neighbourhood $U$ of $x=\mathrm{e}_1$ and 
a constant $C>0$ such that 
\begin{equation}\label{pius2}
   d(x,\exp(t Z)(x)) \leq  C  t^{\frac{1}{w_n-1}} 
            \quad \text{ for all } \exp(tZ)(x) \in U.
\end{equation}

Let us fix a positive parameter $\eps>0$ 
and let $(\gamma^\eps_1, \gamma^\eps_2)$ be the planar curve
obtained by the concatenation of the following three curves: 
the line segment
from $(0,1)$ to $(0, \eps)$,
the curve $(\eps \sigma_1, \eps \sigma_2)$, and 
the line segment from $(\eps , 0)$  to $({1},0)$.
When $a=0$ and $b=1$
we consider the same curve but starting from $(1,0)$.
Let $\gamma^\eps = (\gamma^\eps_1, 
\ldots,   \gamma^\eps_n)$
be the $\mathcal D$-horizontal lift  of this curve to $\R^n$, starting 
from the point $\mathrm{e}_2$ (starting from $\mathrm{e}_1$, when $a=0$ and
$b=1$).
Notice that the $\mathcal D$-horizontal lift 
of $(\eps\sigma_1,\eps\sigma_2)$
is the curve   $\delta_\eps \circ \sigma =
(\eps^{w_1} \sigma_1,\ldots,  \eps^{w_n} \sigma_n )$, 
by   \eqref{dilations2}, 
and hence the  end-point  of $\delta_\eps\circ\sigma$ is
the point $ \eps \mathrm{e_1} +   \eps^{w_n}  h\mathrm{e}_n$.
Moving along  the vector field $X_1$
does not change the
$n$th coordinate $x_n$, hence we conclude that  
the final point of $\gamma^\eps$ is
$x^\eps = \mathrm{e}_1 +  \eps^{w_n}  h \mathrm{e}_n$.
By \eqref{formabj} and \eqref{bubbles} we have  
$\mathrm{e}_1 +  \eps^{w_n}  h \mathrm{e}_n
= \exp(\eps^{w_n} h Z)(\mathrm{e}_1)$. 
Since $x^\eps \to x = \mathrm{e}_1$ as $\eps \to 0$,  
by \eqref{pius2} and for $\eps$ small enough we have  
\begin{equation}\label{D!}
   d  (x,x^\eps)\leq  C
 h ^{\frac{1}{w_n-1}} \eps ^{\frac{w_n}{w_n-1}}.
\end{equation}

The sub-Riemannian length in $(\R^n,\mathcal D ,g)$ of $\gamma^\eps$  is
\begin{equation}\label{LEN}
\begin{split}
L(\gamma^\eps) &= (1-\eps ) L(\gamma	) +  L(\delta_\eps\circ \sigma)\\
&= (1-\eps ) L(\gamma	) +  \eps L(\sigma)
\\
&= L(\gamma) - \eps (L(\gamma)- L(\sigma))
\\
&= L(\gamma) - \eps (L(\widehat \gamma)- L(\widehat \sigma))
\\ 
&= L(\gamma) - \eps k.
\end{split}
\end{equation}

Thus, from \eqref{D!} and \eqref{LEN} we obtain (below we let $y=\mathrm{e}_2$)
\[
\begin{split}
d (y, x)&
\leq d (y, x^\eps)+d (x^\eps, x)\\
&\leq L(\gamma^\eps) + 
C \eps^{\frac{w_n}
{w_n-1}} h^{\frac{1}{w_n-1}}
\\
&= L(\gamma) - \eps k + C \eps^{\frac{w_n}
{w_n-1}} h^{\frac{1}{w_n-1}}.
\end{split}
\]
Since $k>0$, there exists  an $\eps>0$  such that
$C h ^{\frac{1}{w_n-1}} \eps^{\frac{1}{w_n-1}} <k/2$ and hence
\[
d (x, y) <  L(\gamma) - \eps k /2 < L(\gamma).
\]
This proves   that $\gamma$ is not length minimizing in $(\R^n,\mathcal 
D,g)$.
This also concludes the proof by induction of  Theorem  
\ref{thermos2}.

\section{Proof  of Theorem  \ref{main teo}}
\label{S3}

Let  $\Delta =
\mathrm{span}\{  X_1,X_2\}$ be the distribution of planes in $\R^4$
spanned by the vector fields $X_1$ and $X_2$ in \eqref{X_1,X_2}. 
We fix the metric $g$ on $\Delta$ making $X_{1},X_{2}$ an orthonormal frame for $\Delta$. 
Length minimizers for the sub-Riemannian distance are {\em
extremals} in the sense of Geometric Control 
Theory, i.e., they satisfy certain
necessary conditions given by Pontryagin Maximum Principle. 
Extremals may be either normal or abnormal. 
Normal extremals are always smooth.
The following proposition classifies abnormal nonsmooth extremals.

\begin{proposition} \label{abnormal:proposition}
In the structure $(\R^4,\Delta)$,
the only nonsmooth abnormal extremals are the curves 
\begin{equation}\label{horizontal:corner}
\gamma(t) = \begin{cases}
  (-t x_1, -t x_2, 0,a)   & \text{if }t\in [-1,0]
\\ 
(t y_1, t y_2, 0,a)  & \text{if }t\in (0,1],
\end{cases}
\end{equation}
where $a\in \R$ and $(x_1,x_2), (y_1,y_2)\in\R^2$ are linearly independent.
\end{proposition}

\proof[Proof]
Let $\gamma:[0,1]\to\R^4$ be an 
abnormal extremal of the distribution $\Delta$.
By  Pontryagin Maximum Principle, 
there exists an absolutely continuous curve
$\xi:[0,1]\to\R^4$ solving almost everywhere the system of
differential equations
\begin{equation}
\label{equation}
\dot\xi = \big(2 \dot\gamma_2 \xi_3,
\,-2\dot\gamma_1\xi_3, \,-2 \gamma_3 \dot\gamma_1 \xi_4,\, 0\big).
\end{equation}
See, e.g., \cite[Theorem 2.1]{LLMV} for a formulation of
Pontryagin Maximum Principle.
Moreover, we have $\langle X_1(\gamma), \xi\rangle 
=\langle X_2(\gamma),\xi\rangle=\langle [X_1,X_2],\xi\rangle = 0$,
where $\langle \cdot,\cdot\rangle$ is the standard scalar product
of $\R^4$. The last equation $\langle [X_1,X_2],\xi\rangle = 0$
is Goh condition, that holds automatically true in the rank 2 case (see, e.g., \cite{Vittone-note}).
Namely, the curve $\xi$ also solves the system of equations 
\begin{equation}\label{SYS}
 \begin{split}
& \xi_1 + 2 \gamma_2 \xi_3 + \gamma_3^2 \xi_4 = 0 \\
 & \xi_2 - 2 \gamma_1 \xi_3 =0 \\
  & \xi_3-\gamma_1 \gamma_3 \xi_4=0.
 \end{split}
\end{equation}
From \eqref{equation}, we see that  that $\xi_4$ is constant.
This constant is nonzero, otherwise   \eqref{SYS} would
imply $\xi=0$, and this is not possible
for abnormal extremals.
By linearity we can assume that  $\xi_4=1$, and thus
\eqref{SYS} trasforms into the system 
\begin{equation}\label{equation:dual}
\xi = \big ( -2\gamma_1 \gamma_2 \gamma_3-\gamma_3^2, 
  2\gamma_1^2\gamma_3,
\gamma_1\gamma_3, 1\big),
\end{equation}
and the system \eqref{equation} becomes
\begin{equation}
\label{equation1}
\dot\xi = \big(2 \gamma_1\gamma_3 \dot\gamma_2  ,
-2\gamma_1\gamma_3 \dot\gamma_1, -2 \gamma_3 \dot\gamma_1 , 0\big).
\end{equation}
Differentiating the second equation in \eqref{equation:dual}, we find 
$\dot\xi_2 = 4\gamma_1\gamma_3 \dot\gamma_1+2\gamma_1^2\dot\gamma_3$,
and comparing with the second equation in \eqref{equation1}, we deduce that
$\gamma_1^ 2 (3\dot\gamma_1\gamma_3+ \gamma_1 \dot\gamma_3) = 0$.
This in turn implies that the function 
$\phi(t) =  \gamma_1(t)^3 \gamma_1(t)$ is a constant $c\in\R$.

Now there are two cases.

\emph{First case: } $c=0$. In this case, 
the equation $\gamma_1^3\gamma_3=0$ implies that $\gamma$ is either
a line or a corner of the form \eqref{horizontal:corner}.

\emph{Second case: } $c\neq 0$. In this case,
by differentiating the identity $\gamma_3 = c/\gamma_1^3$  and using 
the horizontality condition  
 $\dot\gamma_3=2\dot\gamma_1 \gamma_2-2\dot \gamma_2 \gamma_1$,
we deduce that 
\[
 \big  \langle (3c +2  \gamma_1^4\gamma_2, -2 \dot\gamma_1^5), (\dot\gamma_1,
\dot\gamma_2)\big\rangle=0,
\]
where $\langle\cdot,\cdot\rangle$ is the standard scalar product of $\R^2$.
In other words, the planar curve 
$(\gamma_1, \gamma_2)$ is,
up to reparameterization, an integral curve of the
vector field in the plane 
\[
  \frac{ \partial}{\partial x_1} + \frac{3 c+ 2 x_1^4x_2}{2x_1^5} 
  \frac{ \partial}{\partial x_2},
\quad x_1\neq 0.
\]
Thus the curve $\gamma$ is
 \begin{equation}
 \label{form}
 \gamma(t) =  \Big(t, bt -\frac{3}{10} c t^{-4},
 ct^{-3} 
 ,
-\frac{1}{5} c^2 t^{-5} +d
 \Big), \quad \text{ with } t\neq 0,
 \end{equation}
 for some $b,c,d\in \R$. All such curves are  $C^\infty$. 

We conclude that   the nonsmooth abnormal extremals in $(\R^4,\Delta)$
are precisely the corners \eqref{horizontal:corner}.
  \qed

\begin{remark}
In the proof of Proposition \ref{abnormal:proposition}, 
we have the formula \eqref{equation:dual} for the dual curve $\xi$
of an abnormal extremal $\gamma$. The coordinates of $\xi$ are
polynomial functions of the coordinates of $\gamma$.
This is analogous to the results obtained in 
\cite{LLMV, LLMV2} for stratified nilpotent groups. 
In such groups, 
dual curves can be reconstructed using a special family of 
polynomials, called \emph{extremal polynomials}, and abnormal
extremals are always contained in the level sets of extremal polynomials.
\end{remark}

We conclude with the proof of Theorem \ref{main teo}.
\proof[Proof of Theorem \ref{main teo}]
Thanks to Proposition \ref{abnormal:proposition}, it is enough to prove the non-minimality of corners in $(\R^{4},\Delta)$ at $x=0$. Since the distribution $\Delta$ satisfies the assumption \eqref{LALLA} at $x=0$, we can use Theorem \ref{thermos2} and obtain the desired conclusion.
\qed

\bibliographystyle{amsalpha}

\def\cprime{$'$} \def\cprime{$'$}
\providecommand{\bysame}{\leavevmode\hbox to3em{\hrulefill}\thinspace}
\providecommand{\MR}{\relax\ifhmode\unskip\space\fi MR }
\providecommand{\MRhref}[2]{%
  \href{http://www.ams.org/mathscinet-getitem?mr=#1}{#2}
}
\providecommand{\href}[2]{#2}

\end{document}